\documentclass[reqno]{amsart}
\usepackage{amsfonts,amsmath,amssymb}

\numberwithin{equation}{section}

\newtheorem{thm}{Theorem}[section]

\begin{document}
\title[Stability of the $1$-Laplacian eigenvalue]
{Stability and perturbations of the domain\\
for the first eigenvalue of the $1$-Laplacian}
\author{Emmanuel Hebey}
\address{Emmanuel Hebey, Universit\'e de Cergy-Pontoise, 
D\'epartement de Math\'ematiques, Site de 
Saint-Martin, 2 avenue Adolphe Chauvin, 
95302 Cergy-Pontoise cedex, 
France}
\email{Emmanuel.Hebey@math.u-cergy.fr}
\author{Nicolas Saintier}
\address{Nicolas Saintier, Universit\'e Pierre et Marie Curie, D\'epartement de Math\'ematiques, 
4 place Jussieu, 75252 Paris cedex 05, France}
\email{saintier@math.jussieu.fr}

\date{March 30, 2005}

\maketitle

Let $\Omega$ be a smooth bounded domain in $\mathbb{R}^n$, $n \ge 2$. The $1$-Laplacian 
on $\Omega$ is the formal operator
$$\Delta_1u = -\hbox{div}\left(\frac{\nabla u}{\vert\nabla u\vert}\right)$$
we get by a formal derivation of $F(u) = \int_\Omega\vert\nabla u\vert dx$, or by letting $p \to 1$  
in the definition of the $p$-Laplacian, $p > 1$.  
By analogy with the definition of the first eigenvalue
$\lambda_{p,\Omega}$ of the p-Laplacian on $\Omega$, 
we define the first eigenvalue $\lambda_{1,\Omega}$ of the $1$-laplacian on
$\Omega$ by the minimization problem
\begin{equation}\label{FirstDef}
 \lambda_{1,\Omega} = \inf_{\begin{cases} u\in \dot H_1^1(\Omega) \\
    \int_{\Omega}\vert u\vert dx = 1 \end{cases}} \int_\Omega \vert\nabla u\vert dx\hskip.1cm ,
\end{equation}
where $\dot H_1^1(\Omega)$ is the closure of $C^\infty_0(\Omega)$ in the 
Sobolev space $H_1^1(\Omega)$ of functions in $L^1(\Omega)$ with one derivative in $L^1$. 
By geometric 
measure theory, and the coarea formula, we also have that $\lambda_{1,\Omega} = h(\Omega)$,
where $h(\Omega)$ is defined 
as the infimum of the ratio $\vert\partial D\vert/\vert D\vert$, $D$ varies over all smooth subdomains $D \subset\subset \Omega$, and $\vert\partial D\vert$ and $\vert D\vert$ are the 
$(n-1)$-dimensional and $n$-dimensional measures of $\partial D$ and $D$.  
The result is known as Cheeger's theorem \cite{Che}, and $h(\Omega)$ is known 
as the Cheeger constant of $\Omega$ (see for instance Chavel 
\cite{Cha}). Note the infimum in $h(\Omega)$ 
is not attained by a smooth subdomain $D \subset\subset \Omega$ since, if not, we may
blow it up by a factor larger than one. This would decrease $h$, contradicting the optimality of $D$. 
Minimizers for $h(\Omega)$ touch the boundary $\partial\Omega$.

\medskip The main purpose of this paper is the study of the
dependence of 
$\lambda_{1,\Omega}$ under perturbations of $\Omega$. The notion of
perturbation is here quantified by mean of the $1$-capacity. We provide results
ranging from general type of perturbations to regular perturbations
by diffeomorphisms. This type of problem has been widely studied in the
case of the Laplacian, hardly in the case of the $p$-Laplacian and,
as far as we know, has not been studied before in the
case of the 1-laplacian. 
A natural space to study $\lambda_{1,\Omega}$ is $BV(\Omega)$, the space of 
functions of bounded variations (see, for instance, 
Evans and Gariepi \cite{EvaGar}, or Giusti \cite{Giu}). By standard properties of the space $BV(\Omega)$, 
we can also define $\lambda_{1,\Omega}$ by the equation
\begin{equation}\label{SecondDef}
 \lambda_{1,\Omega} = \inf_{\begin{cases} u\in BV(\Omega) \\
    \int_{\Omega}\vert u\vert dx = 1 \end{cases}} \left(\int_\Omega\vert\nabla u\vert dx 
    + \int_{\partial\Omega}\vert u\vert d\sigma\right)\hskip.1cm .
\end{equation}
Note here that if $u \in BV(\Omega)$, and $\overline{u}$ is the extension of $u$ by $0$ 
in ${\mathbb R}^n\backslash\overline{\Omega}$, then $\overline{u} \in BV({\mathbb R}^n)$ and
\begin{equation}\label{ExtEqt}
\int_{{\mathbb R}^n}\vert\nabla\overline{u}\vert dx = \int_\Omega\vert\nabla u\vert dx 
+ \int_{\partial\Omega}\vert u\vert d\sigma
\end{equation}
By lower semicontinuity of the total variation, and compactness of the 
embedding $BV(\Omega) \subset L^1(\Omega)$, it easily follows from (\ref{ExtEqt}) that 
the infimum in (\ref{SecondDef}) 
is attained by some nonnegative $u \in BV(\Omega)$. Then 
$u$ is a solution of the equation $\Delta_1u = \lambda_{1,\Omega}$ 
in the sense that there exists $\Lambda \in L^\infty(\Omega,{\mathbb R}^n)$, 
$\Vert\Lambda\Vert_\infty \le 1$, such that
\begin{equation}\label{vp11}
  \begin{cases}
 -\hbox{div} \Lambda = \lambda_{1,\Omega} \hskip.1cm ,\hskip.1cm u \ge 0\hskip.1cm ,\\
   \Lambda\nabla u = \vert\nabla u\vert \hskip.1cm\hbox{in}\hskip.1cm
   \Omega\hskip.1cm ,\hskip.1cm\hbox{and} \\
   (\Lambda\nu)u = -\vert u\vert \hskip.1cm\hbox{on}\hskip.1cm \partial\Omega\hskip.1cm ,
  \end{cases}
 \end{equation}
where $\nu$ is the unit outer normal to $\partial\Omega$, and $ \Lambda\nabla u $ 
is the distribution defined by integrating by parts $\int_\Omega(\Lambda\nabla u)vdx$ when 
$v \in C^\infty_0(\Omega)$ and $\hbox{div}\Lambda$ makes sense. Moreover, see 
for instance Demengel \cite{Dem}, 
$u \in L^\infty(\Omega)$ and, since $\Vert u\Vert_1 = 1$, there exists $C = C(n,\Omega)$, $C > 0$, such that 
$\Vert u\Vert_\infty \le C$. We say $u$ is an eigenfunction for $\lambda_{1,\Omega}$. Now 
define a Cacciopoli set in $\Omega$ as a set $D \subset \Omega$ such that $\chi_D \in BV(\Omega)$, where $\chi_D$ is the characteristic function of $D$. 
Since $\lambda_{1,\Omega} = h(\Omega)$, there are 
Caccioppoli sets $D \subset \Omega$, e.g the level sets of eigenfunctions,
such that $u = \vert D\vert^{-1}\chi_D$ is a minimizer for the right hand side in  
(\ref{SecondDef}).  Such sets are referred to  
as eigensets for $\lambda_{1,\Omega}$ (and sometimes also as Cheeger's sets).  
A general discussion about uniqueness and nonuniqueness of eigensets is in Fridman and Kawohl \cite{FriKaw}. 
We refer also to Stredulinsky and Ziemer \cite{StrZie} (and to Belloni and Kawohl \cite{BelKaw} for the $p$-Laplace 
case). 
Concerning regularity, possible references are Almgren \cite{Alm}, 
De Giorgi \cite{Gio}, 
Gonzales, Massari and Tamanini \cite{GonMasTam}, and Stredulinsky 
and Ziemer \cite{StrZie}. By symmetrization, see Fridman and Kawohl \cite{FriKaw} for details,
$$\lambda_{1,\Omega} 
\ge n\omega_n^{1/n}\vert\Omega\vert^{-1/n}\hskip.1cm ,$$
where $\omega_n$ is the volume of the unit $n$-sphere.

\medskip Let $K$ be a compact subset of ${\mathbb R}^n$. The $1$-capacity of $K$, 
denoted by $\hbox{cap}_1(K)$, is defined as the infimum of the $L^1$-norm of $\vert\nabla u\vert$, 
where the infimum is taken over all $u \in C^\infty_0({\mathbb R}^n)$ such that 
$K \subset \hbox{int}\left\{u \ge 1\right\}$. 
Another possible definition (see, for instance, Maz'ja \cite{Maz}) is that $\hbox{cap}_1(K) = \inf\vert\partial\omega\vert$, 
where the infimum is taken over all smooth open bounded subset $\omega$
such that $K\subset\omega$. In particular, by the isoperimetric inequality, $\vert K\vert^{(n-1)/n} 
\le C\hbox{cap}_1(K)$, where $C > 0$ does not depend on $K$ (but only on the dimension). 
For $A$ and 
$B$ two subsets of ${\mathbb R}^n$, we denote by $A\Delta B$ the symmetric difference of $A$
and $B$. Namely, 
$A\Delta B = (A\backslash B)\cup(B\backslash A)$. 
Our first 
result deals with general type of perturbations of a domain. It states as follows.

\begin{thm}\label{TH1}
Let $\Omega$ be a smooth bounded
 domain in ${\mathbb R}^n$, 
 $(\Omega_\delta)_{\delta > 0}$ be a sequence of smooth bounded domains in $\mathbb{R}^n$, 
and $K_\delta = \hbox{adh}(\Omega\Delta\Omega_\delta)$ be the closure of the symmetric difference $\Omega\Delta\Omega_\delta$. Let $(A_\delta)$ be a sequence of eigensets for the 
$\lambda_{1,\Omega_\delta}$'s. Assume 
$\hbox{cap}_1(K_\delta) 
\to 0$ as $\delta \to 0$. Then, for any $\delta$,
\begin{equation}\label{EqtTh1}
\left\vert\lambda_{1,\Omega_\delta} - \lambda_{1,\Omega}\right\vert = 
\frac{\varepsilon_\delta}{\vert A_\delta\vert}\hbox{cap}_1(K_\delta) + o\left(\hbox{cap}_1(K_\delta)\right)\hskip.1cm ,
\end{equation}
where $\varepsilon_\delta \in [0,1]$ for all $\delta$, 
and $\vert A_\delta\vert \ge \Lambda$ for some $\Lambda > 0$ and all $\delta$. 
In particular, $\lambda_{1,\Omega_\delta} \to \lambda_{1,\Omega}$ 
as $\delta \to 0$. Moreover, up to a subsequence, 
$\chi_{A_\delta} \to \chi_A$ in $L^1({\mathbb R}^n)$ as $\delta \to 0$, where $A$ is an eigenset  
for $\lambda_{1,\Omega}$.
\end{thm}

We stated the second part of Theorem \ref{TH1} for characteristic functions of eigensets. 
However, note the convergence holds also for eigenfunctions. 
In Section \ref{ProofTH1} we prove that if the $u_\delta$'s are 
eigenfunctions for $\lambda_{1,\Omega_\delta}$, then, up to a subsequence, 
with the notation in  (\ref{ExtEqt}), 
$\overline{u}_\delta \to \overline{u}$ in $L^1({\mathbb R}^n)$ as $\delta \to 0$, where $u$ is an eigenfunction 
for $\lambda_{1,\Omega}$. We also get that the measures $\mu_\delta = 
\vert\nabla\overline{u}_\delta\vert$ and $\mu = \vert\nabla\overline{u}\vert$ satisfy $\mu_\delta \rightharpoonup \mu$ weakly 
as $\delta \to 0$.

\medskip A particular case of the general perturbations considered in Theorem \ref{TH1} 
is when we consider domains with holes.  Such domains (in the case of one ball) 
were considered by Sango \cite{San} when 
discussing the first eigenvalue of the $p$-Laplace operator, $p > 1$. 
For $A$ a 
Caccioppoli set, we let $\partial^\star A$ be its reduced boundary, namely (see \cite{EvaGar} or \cite{Giu}) 
the subset of the boundary which is $C^1$ in a measure
theoretic sense. For $x \in {\mathbb R}^n$ and $r > 0$, we let also $B_x(r)$ be the 
$n$-dimensional ball of center $x$ and radius $r$, and $b_n$ be the volume of the unit $n$-dimensional ball.

\begin{thm}\label{TH2}
Let $\Omega\subset\mathbb{R}^n$ be a smooth bounded domain in ${\mathbb R}^n$, 
$A$ be an eigenset for $\lambda_{1,\Omega}$, $x_1,\dots, x_k$ be points in $\Omega$, 
$(\varepsilon_{i,\delta})_{\delta > 0}$, $i = 1,\dots,k$,  be $k$ sequences of positive real numbers 
converging to $0$ as $\delta \to 0$, 
 $K_\delta = \bigcup_{i=1}^k\overline{B}_{x_i}(\varepsilon_{i,\delta})$, $\delta > 0$ small, and 
 $\Omega_\delta = \Omega\backslash K_\delta$. 
 Assume there exists $i_0 = 1,\dots,k$ such that 
 $\varepsilon_{i,\delta} = o\left(\varepsilon_{i_0,\delta}\right)$ for 
 all $i = 1,\dots,k$, $i \not= i_0$. Then, 
 $\hbox{cap}_1(K_\delta) = \omega_{n-1}\varepsilon_{i_0,\delta}^{n-1} 
 + o(\varepsilon_{i_0,\delta}^{n-1})$, and
 \begin{equation}\label{TH2Eqt1}
\lambda_{1,\Omega} \le \lambda_{1,\Omega_\delta} \le 
\lambda_{1,\Omega} + \frac{\omega_{n-1}}{\vert A\vert}\varepsilon_{i_0,\delta}^{n-1}  + 
o\left(\varepsilon_{i_0,\delta}^{n-1}\right)
\end{equation}
for all $\delta$. Moreover,
 \begin{equation}\label{TH2Eqt2}
 \lambda_{1,\Omega_\delta} \le \lambda_{1,\Omega} +
   \begin{cases}
     o\left(\varepsilon_{i_0,\delta}^{n-1}\right) \hskip.1cm\hbox{if}\hskip.1cm x_{i_0}\in \hbox{int}(\Omega\backslash A)\\
    \frac{\omega_{n-1}-2b_{n-1}}{2\vert A\vert}\varepsilon_{i_0,\delta}^{n-1} +
    o\left(\varepsilon_{i_0,\delta}^{n-1}\right) \hskip.1cm\hbox{if}\hskip.1cm  x_{i_0}\in \partial^\star A\hskip.1cm , \\
    \end{cases}
 \end{equation}
where $\hbox{int}(\Omega\backslash A)$ is the interior of $\Omega\backslash A$, and 
$\partial^\star A$ is the reduced boundary of A. 
\end{thm}

Needless to say, 
$\hbox{cap}_1(K_\delta) \to 0$ as $\delta \to 0$ and the convergence of eigensets (resp. eigenfunctions) stated in Theorem \ref{TH1} holds true. Characterizations in dimension $2$ of  
convex $\Omega$'s for which $\Omega = A$ 
(as a by-product for which $\Omega\backslash A \not= \emptyset$, and thus for which equation (\ref{TH2Eqt2}) is not empty) are 
in Bellettini, Caselles and Novaga \cite{BelCasNov}, and Kawohl and Lachand-Robert \cite{KawLac}. When $\Omega$ is 
convex and $n = 2$, the eigenset $A = A_\Omega$ is unique. The domain $\Omega$ is 
said to be a calibrable set when $A = \Omega$. The ovoid domain 
$(x^2+y^2)^2 < x^3$ is a nice example in Kawohl and Lachand-Robert \cite{KawLac} of a 
noncalibrable set. With respect to the notation in Theorem \ref{TH1}, (\ref{TH2Eqt2}) gives that 
we can take $\varepsilon_\delta = 0$ if $x_{i_0}\in \hbox{int}(\Omega\backslash A)$, and $\varepsilon_\delta < 1/2$ if 
$x_{i_0} \in \partial^\star A$. On the other hand, 
when $x_{i_0} \in \hbox{int}(A)$, the upper bound in (\ref{TH2Eqt1}), 
where $\varepsilon_\delta = 1$, cannot (in general) be improved. If we let 
$\Omega = B_0(r)$ and $K_\varepsilon = \overline{B}_0(\varepsilon)$ 
be $n$-dimensional balls, $r > 0$, $0 < \varepsilon \ll 1$, then both 
$\Omega$ and the annulus $\Omega_\varepsilon = \Omega\backslash K_\varepsilon$ are calibrable sets 
(see, for instance, 
Kawohl and Lachand-Robert \cite{KawLac}, Demengel, De Vuyst, and Motron \cite{DemDevMot}). 
In particular, the eigenvalues 
$\lambda_{1,\Omega}$ and $\lambda_{1,\Omega_\varepsilon}$ are given by
$\lambda_{1,\Omega} = n/r$ and $\lambda_{1,\Omega_\varepsilon} = n(r^{n-1}+\varepsilon^{n-1})/(r^n-\varepsilon^n)$, and 
we get that
$$\lambda_{1,\Omega_\varepsilon} = \lambda_{1,\Omega} + \frac{\omega_{n-1}}{\vert A\vert} \varepsilon^{n-1} 
+ o\left(\varepsilon^{n-1}\right)$$
for all $\varepsilon$, where $A = \Omega$ is the eigenset of $\lambda_{1,\Omega}$. Without further 
assumptions, the upper bound in (\ref{TH2Eqt1}) is sharp.

\medskip We now consider the case of a regular perturbation of $\Omega$ by
diffeomorphisms. We prove the differentiability at 0 of the map
$\delta\to\lambda_{1,\Omega_\delta}$ and the convergence of the
eigenfunctions without assumptions on the capacity of
$K_\delta$. In the case of the $p$-Laplacian, $p > 1$, such type of problems have been considered 
by Lamberti \cite{Lam} and Garcia Melian and Sabina
De Lis \cite{GarSab}. For $\Omega$ a smooth bounded open subset of 
${\mathbb R}^n$, we let $(T_\delta)_\delta$ be a family of $C^1$-diffeomorphisms of the form
\begin{equation}\label{DiffeomDefEqt}
T_\delta(x) = \left(1 - \delta\Lambda\right) x + R(x,\delta)\hskip.1cm ,
\end{equation}
where $x \in \overline{\Omega}$, $\Lambda \in {\mathbb R}$, $\delta \in \left(-\delta_0,\delta_0\right)$, $\delta_0 > 0$ is small, 
and $R(.,\delta) \in C^1(\overline{\Omega},\mathbb{R}^n)$ is a perturbative term such that 
 $R(x,\delta) = o(\delta)$ and $D_xR(x,\delta) = o(\delta)$ as $\delta \to 0$, uniformly in $x$. 
In particular, $R(x,0) = 0$, and if $\Omega_\delta = T_\delta(\Omega)$, then 
$\Omega_0 = \Omega$.
 
\begin{thm} \label{TH3}
 Let $\Omega$ be a smooth bounded open subset of $\mathbb{R}^n$, and
 $\Omega_\delta = T_\delta(\Omega)$, where the $T_\delta$'s are
 $C^1$-diffeomorphisms like in (\ref{DiffeomDefEqt}). The 
 function $\delta \to \lambda_{1,\Omega_\delta}$ is continuous and differentiable at 
 $\delta = 0$, and 
$\left(\lambda_{1,\Omega_\delta}\right)^\prime(0) = \Lambda \lambda_{1,\Omega}$, 
 where $\Lambda$ is as in (\ref{DiffeomDefEqt}).
\end{thm}

In these examples of Theorem \ref{TH3}, the $1$-capacity of $K_\delta = \hbox{adh}(\Omega\Delta\Omega_\delta)$ 
can be large. For instance, if $\Omega = B_0(r)$ and $T_\delta = (1-\delta\Lambda)x$, $\Lambda \not= 0$, 
then $K_\delta$ is an annulus with inner (or outer, depending on the sign of 
$\Lambda$) radius $r$. In particular, 
$\hbox{cap}_1(K_\delta) = \omega_{n-1}r^{n-1} + o(1)$, 
and $\hbox{cap}_1(K_\delta) \not\to 0$ as $\delta \to 0$. In other cases, the $1$-capacity of $K_\delta$ 
may tend to zero, and we are back to the situation studied in Theorem \ref{TH1}. For 
instance, if $T_\delta(x) = \left(1 + \delta^\alpha R(\frac{1}{\delta}x)\right) x$, $\alpha > 2$, 
$R \in C^1_0({\mathbb R}^n)$, and $0 \in \partial\Omega$, then $K_\delta \subset B_0(r_0\delta)$ 
for some $r_0 > 0$, and $\hbox{cap}_1(K_\delta) \to 0$ as $\delta \to 0$. Combining Theorems \ref{TH1} 
and \ref{TH3}, $\Lambda = 0$ if $\hbox{cap}_1(K_\delta) = o(\delta)$.
The following sections are devoted to the proofs of the above theorems.

\section{proof of theorem \ref{TH1}}\label{ProofTH1}

Let $A_\delta$ be an eigenset of $\lambda_{1,\Omega_\delta}$, $\delta > 0$ fixed. 
Let also $\omega_\delta$ be a smooth 
bounded open subset such that $K_\delta \subset \omega_\delta$ and
\begin{equation}\label{ProofTH1Eqt2}
\hbox{cap}_1(K_\delta) \le \vert\partial\omega_\delta\vert \le 
\hbox{cap}_1(K_\delta) + \varepsilon_\delta\hskip.1cm ,
\end{equation}
where $\varepsilon_\delta > 0$ is such that $\varepsilon_\delta = o\left(\hbox{cap}_1(K_\delta)\right)$. 
We define $v_\delta$ by $v_\delta = \chi_{A_\delta}$ in ${\mathbb R}^n\backslash\overline{\omega}_\delta$, 
and $v_\delta = 0$ in $\omega_\delta$. Then $\hbox{supp}v_\delta \subset \overline{\Omega}$, where 
$\hbox{supp}v_\delta$ is the support of $v_\delta$. Since $v_\delta \le 1$, 
we can write that
\begin{equation}\label{ProofTH1Eqt3}
\begin{split}
\int_{{\mathbb R}^n}\vert\nabla v_\delta\vert dx
&\le \int_{{\mathbb R}^n}\vert\nabla\chi_{A_\delta}\vert dx 
+ \vert\partial\omega_\delta\vert\\
&=  \vert A_\delta\vert\lambda_{1,\Omega_\delta}
 + \left(1+o(1)\right)\hbox{cap}_1(K_\delta)\hskip.1cm ,
\end{split}
\end{equation}
where $o(1) \to 0$ as $\delta \to 0$. We can also write that
\begin{equation}\label{ProofTH1Eqt6}
\begin{split}
\int_\Omega v_\delta dx 
& = \vert A_\delta\vert - \int_{\omega_\delta}\chi_{A_\delta}dx\\
& = \vert A_\delta\vert + O\left(\vert\omega_\delta\vert\right)\hskip.1cm .
\end{split}
\end{equation}
By the isoperimetric inequality, and by  (\ref{ProofTH1Eqt2}),
$\vert\omega_\delta\vert
\le C \vert\partial\omega_\delta\vert^{\frac{n}{n-1}}
= o\left(\hbox{cap}_1(K_\delta)\right)$, where 
$C > 0$ is a dimensional constant independent of $\delta$. Coming back to 
(\ref{ProofTH1Eqt6}), it follows that
\begin{equation}\label{ProofTH1Eqt9}
\int_\Omega v_\delta dx = \vert A_\delta\vert 
+ o\left(\hbox{cap}_1(K_\delta)\right)
\end{equation}
and by the variational definition of $\lambda_{1,\Omega}$, 
 (\ref{ProofTH1Eqt3}), and (\ref{ProofTH1Eqt9}), we get that
\begin{equation}\label{ProofTH1Eqt10}
\lambda_{1,\Omega} \le \lambda_{1,\Omega_\delta} + \frac{1}{\vert A_\delta\vert} \hbox{cap}_1(K_\delta) 
+ o\left(\hbox{cap}_1(K_\delta)\right)
\end{equation}
for all $\delta > 0$. Similar arguments give that the converse inequality holds also. 
We let $A$ be an eigenset for $\lambda_{1,\Omega}$, and let 
$w_\delta$ be the function given by $w_\delta = \chi_A$ in ${\mathbb R}^n\backslash\overline{\omega}_\delta$, 
and $w_\delta = 0$ in $\omega_\delta$, where $\omega_\delta$ is as in 
(\ref{ProofTH1Eqt2}). Then $\hbox{supp}w_\delta \subset \overline{\Omega}_\delta$, and, as above, we can write
\begin{equation}\label{ProofTH1Eqt15}
\begin{split}
&\int_{{\mathbb R}^n}\vert\nabla w_\delta\vert dx
\le  \vert A\vert\lambda_{1,\Omega} + \left(1+o(1)\right) \hbox{cap}_1(K_\delta)
\hskip.1cm ,\hskip.1cm \hbox{and}\\
&\int_{\Omega_\delta}w_\delta dx = \vert A\vert + o\left(\hbox{cap}_1(K_\delta)\right)\hskip.1cm .
\end{split}
\end{equation}
In particular, it follows from the variational definition of $\lambda_{1,\Omega_\delta}$ that
\begin{equation}\label{ProofTH1Eqt16}
\lambda_{1,\Omega_\delta} \le \lambda_{1,\Omega} + \frac{1}{\vert A\vert} \hbox{cap}_1(K_\delta) 
+ o\left(\hbox{cap}_1(K_\delta)\right)
\end{equation}
for all $\delta > 0$. 
Without loss of generality, by the lower semicontinuity of the total variation, we may choose
$A = A_0$ in (\ref{ProofTH1Eqt16}) such that it is of maximum volume among the eigensets for $\lambda_{1,\Omega}$.

\medskip In what follows we let $u_\delta \in BV(\Omega_\delta)$ be an eigenfunction for $\lambda_{1,\Omega_\delta}$, 
like for instance $u_\delta = \vert A_\delta\vert^{-1}\chi_{A_\delta}$, 
and we assume that $\hbox{cap}_1(K_\delta) \to 0$ as $\delta \to 0$. 
Then, by (\ref{ProofTH1Eqt10}) and (\ref{ProofTH1Eqt16}), 
$\lambda_{1,\Omega_\delta} \to \lambda_{1,\Omega}$ as $\delta \to 0$. 
Note here that 
$\vert A_\delta\vert \ge C$ for some $C > 0$ (thanks for instance to the 
Sobolev inequality in $BV({\mathbb R}^n)$ that we apply to the extensions 
by zero outside $\overline{\Omega}_\delta$ of the functions $\vert A_\delta\vert^{-1}\chi_{A_\delta}$). 
In what follows we let $D$ be a smooth bounded open subset of 
${\mathbb R}^n$ such that $\overline{\Omega} \subset D$, and 
$\overline{\Omega}_\delta \subset D$ for all $\delta$. We let $\overline{u}_\delta$ 
be the extension of $u_\delta$ by zero outside $\overline{\Omega}_\delta$. 
By (\ref{ProofTH1Eqt16}), the sequence 
$(\overline{u}_\delta)$ is bounded in $BV(D)$. Then, by compactness of the embedding 
of $BV(D)$ into $L^1(D)$, we may assume that, up to a subsequence, $\overline{u}_\delta \to 
\overline{u}$ in $L^1(D)$ for some $\overline{u} \in BV(D)$. By the Sobolev inequality for $BV$-functions, we also have that the $\overline{u}_\delta$'s are bounded in the 
Lebesgue's space $L^{n/(n-1)}(D)$. 
On the one hand, we have that
\begin{eqnarray*}
\int_{D\backslash\Omega}\vert\overline{u}_\delta\vert dx
& \le & \Vert\overline{u}_\delta\Vert_{L^{n/(n-1)}}\vert\Omega_\delta\backslash\Omega\vert^{1/n}\\
& \le & C \vert\Omega_\delta\backslash\Omega\vert^{1/n}\hskip.1cm .
\end{eqnarray*}
On the other hand, we can write with the isoperimetric inequality that
$$\vert\Omega_\delta\backslash\Omega\vert \le \vert K_\delta\vert \le C\hbox{cap}_1(K_\delta)^{\frac{n}{n-1}} 
= o\left(\hbox{cap}_1(K_\delta)\right)\hskip.1cm ,$$
where, as above, $C > 0$ is independent of $\delta$. In particular,  
$\int_{D\backslash\Omega}\vert\overline{u}_\delta\vert dx = o(1)$. We regard $\overline{u}$ as a function in ${\mathbb R}^n$ 
(by letting $\overline{u} = 0$ outside $D$), and let $u = \overline{u}_{\vert\Omega}$ be the restriction 
of $\overline{u}$ to $\Omega$. Then, according to what we just said, 
$\overline{u} = u$ in $\Omega$, and $\overline{u} = 0$ in 
${\mathbb R}^n\backslash\overline{\Omega}$.  As is easily checked, $\int_\Omega\vert u\vert dx = 1$, 
while by lower semicontinuity of the total variation, 
and since $\lambda_{1,\Omega_\delta} \to \lambda_{1,\Omega}$, we can write that
$$\lambda_{1,\Omega} 
= \int_{{\mathbb R}^n}\vert\nabla\overline{u}_\delta\vert dx + o(1)\\
\ge \int_{{\mathbb R}^n}\vert\nabla\overline{u}\vert dx + o(1)
\hskip.1cm .$$
In particular, $u$ is an eigenfunction for $\lambda_{1,\Omega}$, and
\begin{equation}\label{ProofTH1ConvMeasureEst}
\lim_{\delta \to 0}\int_{{\mathbb R}^n}\vert\nabla\overline{u}_\delta\vert dx = 
\int_{{\mathbb R}^n}\vert\nabla\overline{u}\vert dx\hskip.1cm .
\end{equation}
Letting $u_\delta = \vert A_\delta\vert^{-1}\chi_{A_\delta}$, we may assume $\vert A_\delta\vert 
\to \Lambda$ for some $\Lambda > 0$, and $\overline{u}_\delta \to \overline{u}$ a.e. 
Since $\int_\Omega\vert u\vert dx = 1$, we can write that 
$\overline{u} = \vert A\vert^{-1}\chi_A$ for some $A \subset \overline{\Omega}$. 
In particular,  $A$ is an eigenset for $\lambda_{1,\Omega}$, and $\vert A_\delta\vert \to 
\vert A\vert$ so that, by (\ref{ProofTH1Eqt10}) and (\ref{ProofTH1Eqt16}),
$$\lambda_{1,\Omega_\delta} = \lambda_{1,\Omega} + \frac{\varepsilon_\delta}{\vert A_\delta\vert}
\hbox{cap}_1(K_\delta) + o\left(\hbox{cap}_1(K_\delta)\right)\hskip.1cm ,$$
where $\varepsilon_\delta \in [-1,1]$ (since $\vert A\vert \ge \vert A_0\vert$). 
This is equation (\ref{EqtTh1}) in Theorem \ref{TH1}, and the equation holds for all $\delta$ by contradiction. 
Noting that the convergence 
$\overline{u}_\delta \to \overline{u}$ in $L^1$ gives that $\chi_{A_\delta} \to \chi_A$ in 
$L^1$, Theorem \ref{TH1} is proved.

\medskip For the remark following Theorem \ref{TH1}, there is still 
to prove that if $\mu_\delta = \vert\nabla\overline{u}_\delta\vert$ 
and $\mu = \vert\nabla\overline{u}\vert$, then $\mu_\delta \rightharpoonup \mu$ weakly as $\delta \to 0$. 
By lower semicontinuity of the total variation,
$\mu(U) \le \liminf_{\delta \to 0}\mu_\delta(U)$ 
for all open subset $U$ of ${\mathbb R}^n$. Conversely, let us assume that there exists a compact subset 
$K$ of ${\mathbb R}^n$, and $\varepsilon > 0$ such that $\mu(K) + \varepsilon \le \mu_\delta(K)$ for a 
subsequence of the $\mu_\delta$'s. Let $\Omega^\prime$ be an open subset of ${\mathbb R}^n$ which contain 
$\overline{\Omega}$, $K$, and the $\overline{\Omega}_\delta$'s. Up to passing to another 
subsequence, by lower semicontinuity of the total variation, we can assume that 
$\mu(\Omega^\prime\backslash K) \le \mu_\delta(\Omega^\prime\backslash K) + \varepsilon^\prime$ for 
all $\delta$, where $\varepsilon^\prime < \varepsilon$ is positive. Then, if 
$\hat\varepsilon = \varepsilon-\varepsilon^\prime$, we can write that
\begin{eqnarray*}
\mu({\mathbb R}^n) 
= \mu(\Omega^\prime)
& =  & \mu(K) + \mu(\Omega^\prime\backslash K)\\
& \le & \mu_\delta(K) - \varepsilon + \mu_\delta(\Omega^\prime\backslash K) + \varepsilon^\prime\\
& = & \mu_\delta(\Omega^\prime) - \hat\varepsilon
= \mu_\delta({\mathbb R}^n) - \hat\varepsilon
\end{eqnarray*}
for all $\delta$, and we get a contradiction with (\ref{ProofTH1ConvMeasureEst}) 
since $\hat\varepsilon > 0$. As a consequence, for any
compact subset $K$ of ${\mathbb R}^n$,
$\mu(K) \ge \limsup_{\delta \to 0}\mu_\delta(K)$
and, see for instance Evans-Gariepy \cite{EvaGar}, we actually proved that the measures $\mu_\delta$ converge 
weakly to the measure $\mu$. 

\section{Proof of Theorem \ref{TH2}}\label{ProofTH2}

We now turn our attention to the proof of theorem
\ref{TH2}. As is easily checked from the definition of the $1$-capacity, the fact that the $1$-capacity is an outer measure, 
and the isoperimetric inequality in Euclidean space,
\begin{equation}\label{CapEqtProofTH2}
\hbox{cap}_1(K_\delta) = \omega_{n-1}\varepsilon_{i_0,\delta}^{n-1} 
 + o\left(\varepsilon_{i_0,\delta}^{n-1}\right)\hskip.1cm .
 \end{equation}
 From independent considerations, we clearly have that $\dot H_1^1(\Omega_\delta) 
 \subset \dot H_1^1(\Omega)$. Hence, $\lambda_{1,\Omega} \le \lambda_{1,\Omega_\delta}$. 
 On the other hand,  by the proof of Theorem \ref{TH1}, see 
 (\ref{ProofTH1Eqt16}), and by (\ref{CapEqtProofTH2}), 
 we also have that
 \begin{eqnarray*}
 \lambda_{1,\Omega_\delta} 
 & \le & \lambda_{1,\Omega} 
 + \frac{1}{\vert A\vert} \hbox{cap}_1(K_\delta) + o\left(\hbox{cap}_1(K_\delta)\right)\\
 & = & \lambda_{1,\Omega} + \frac{\omega_{n-1}}{\vert A\vert}\varepsilon_{i_0,\delta}^{n-1}  + 
o\left(\varepsilon_{i_0,\delta}^{n-1}\right)\hskip.1cm .
 \end{eqnarray*}
 This proves (\ref{TH2Eqt1}). It remains to prove (\ref{TH2Eqt2}). For this, we need to be more careful 
 than in the proof of Theorem \ref{TH1}. 
 We let $A$ be an eigenset for $\lambda_{1,\Omega}$, and 
 let $\omega_\delta$ be the union from $i=1$ to $k$ of balls $B_{x_i}(\tilde\varepsilon_{i,\delta})$, 
 $\delta > 0$ small, where $\varepsilon_{i,\delta} < \tilde\varepsilon_{i,\delta}$, 
 and $\tilde\varepsilon_{i,\delta} = \left(1+o(1)\right)\varepsilon_{i,\delta}$ for all $i$ and $\delta$. 
 For $u = \chi_A$, 
 we let also $u_\delta^+$ be the trace of $u$ when $u$ is restricted to ${\mathbb R}^n\backslash\overline{\omega}_\delta$, 
 and $u_\delta^-$ be the trace of $u$ when $u$ is restricted to $\omega_\delta$. Then, $\vert A\vert\lambda_{1,\Omega}
 = \int_{{\mathbb R}^n}\vert\nabla u\vert dx$, and 
 \begin{equation}\label{ProofTH2RefinEqt1}
 \int_{{\mathbb R}^n}\vert\nabla u\vert dx
 = \int_{{\mathbb R}^n\backslash\overline{\omega}_\delta}\vert\nabla u\vert dx 
 + \int_{\omega_\delta}\vert\nabla u\vert dx + \int_{\partial\omega_\delta}\vert u_\delta^+-u_\delta^-\vert d\sigma\hskip.1cm .
 \end{equation}
In particular, if we let $w_\delta$ be given by $w_\delta = u$ in ${\mathbb R}^n\backslash\overline{\omega}_\delta$, 
and $w_\delta = 0$ in $\omega_\delta$, we get with (\ref{ProofTH2RefinEqt1}) that 
\begin{equation}\label{ProofTH2RefinKeyEqt}
\begin{split}
\int_{{\mathbb R}^n}\vert\nabla w_\delta\vert dx
&= \int_{{\mathbb R}^n\backslash\overline{\omega}_\delta}\vert\nabla u\vert dx 
+ \int_{\partial\omega_\delta}u_\delta^+d\sigma\\
&\le \vert A\vert\lambda_{1,\Omega} - \int_{\omega_\delta}\vert\nabla u\vert dx 
+ \int_{\partial\omega_\delta}u_\delta^-d\sigma\\
&\le \vert A\vert\lambda_{1,\Omega} - \int_{B_\delta}\vert\nabla u\vert dx 
+ \int_{\partial B_\delta}u_\delta^-d\sigma + o\left(\varepsilon_{i_0,\delta}^{n-1}\right)\hskip.1cm ,
\end{split}
\end{equation} 
where $B_\delta = B_{x_{i_0}}(\tilde\varepsilon_{i_0,\delta})$. We also have that
\begin{equation}\label{ProofTH2RefinKeyEqtContin}
\begin{split}
&\int_{\partial B_\delta}u_\delta^-d\sigma = \vert\partial B_\delta\vert - \int_{\partial B_\delta}(1-u_\delta^-)d\sigma
\hskip.1cm ,\hskip.1cm\hbox{and}\\
&\int_{B_\delta}\vert\nabla u\vert dx + \int_{\partial B_\delta}(1-u_\delta^-)d\sigma 
= \int_{{\mathbb R}^n}\vert\nabla v\vert dx\hskip.1cm ,
\end{split}
\end{equation}
where $v$ is the function $v = \chi_{A^c}$ in $B_\delta$, $v = 0$ in ${\mathbb R}^n\backslash\overline{B}_\delta$, 
and $A^c = {\mathbb R}^n\backslash A$. If we assume that $x_{i_0} \in \hbox{int}(\Omega\backslash A)$, then $u_\delta^- = 0$ 
on $\partial B_\delta$, and it follows from the second equation in (\ref{ProofTH1Eqt15}) of Section 
\ref{ProofTH1}, from (\ref{CapEqtProofTH2}) and (\ref{ProofTH2RefinKeyEqt}), 
and from the variational definition of $\lambda_{1,\Omega_\delta}$, 
that the first equation in (\ref{TH2Eqt2}) is true. Now we assume that $x_{i_0} \in \partial^\star A$. Then, 
the second equation in (\ref{ProofTH1Eqt15}) of Section 
\ref{ProofTH1}, (\ref{CapEqtProofTH2}) , (\ref{ProofTH2RefinKeyEqt})--(\ref{ProofTH2RefinKeyEqtContin}), 
and the variational definition of $\lambda_{1,\Omega_\delta}$ give that
\begin{equation}\label{ProofTH2ContInterEqt}
\lambda_{1,\Omega_\delta} \le \lambda_{1,\Omega} 
+ \frac{1}{\vert A\vert}\left(\omega_{n-1}\varepsilon_{i_0,\delta}^{n-1} - \int_{{\mathbb R}^n}\vert\nabla v\vert dx\right) 
+ o\left(\varepsilon_{i_0,\delta}^{n-1}\right)
\end{equation}
for all $\delta$. Let $T_\delta$ be the diffeomorphism given by $T_\delta(x) = x_{i_0} + \tilde\varepsilon_{i_0,\delta}^{-1}(x-x_{i_0})$. 
Then, by the change of variables formula for the total variation, see Giusti \cite{Giu} 
or equation (\ref{ChangeVarFor}) below, we can write that
\begin{equation}\label{ChgeVarForProofTH2}
\int_{{\mathbb R}^n}\vert\nabla v\vert dx = 
\tilde\varepsilon_{i_0,\delta}^{n-1}
 \int_{{\mathbb R}^n}\vert\nabla(v\circ T_\delta^{-1})\vert dx\hskip.1cm .
 \end{equation}
 Let $A_\delta$ be the set consisting of the $x$ such that $T_\delta^{-1}(x) \in A^c$. Then 
 $v\circ T_\delta^{-1} = \chi_{A_\delta}$ in $B$, and $v\circ T_\delta^{-1} = 0$ in ${\mathbb R}^n\backslash\overline{B}$, 
 where $B = B_{x_{i_0}}(1)$. Since 
 $\partial^\star A^c = \partial^\star A$, $x_{i_0} \in \partial^\star A^c$. By the blow-up property of 
 the reduced boundary (see, for instance, Evans-Gariepy \cite{EvaGar}), we can write that
\begin{equation}\label{BlowUpProofTH2}
\chi_{A_\delta} \to \chi_{H^-(x_{i_0})}\hskip.2cm\hbox{in}\hskip.1cm L^1_{loc}({\mathbb R}^n)
\end{equation}
as $\delta \to 0$, where $H^-(x_{i_0})$ consists of the $y \in {\mathbb R}^n$ 
such that $\nu_{A^c}(x_{i_0}).(y-x_{i_0}) \le 0$, 
and where $\nu_{A^c}(x_{i_0})$ is the 
generalized exterior normal to $A^c$ at $x_{i_0}$. By (\ref{BlowUpProofTH2}), 
$v\circ T_\delta^{-1} \to \hat v$ in $L^1({\mathbb R}^n)$, where $\hat v = \chi_{H^-(x_{i_0})}$ in $B$, 
and $\hat v = 0$ in ${\mathbb R}^n\backslash\overline{B}$. By lower semicontinuity of the total variation, it follows that
\begin{equation}\label{EqtAfterBUProofTH2}
 \int_{{\mathbb R}^n}\vert\nabla(v\circ T_\delta^{-1})\vert dx 
 \ge \int_{{\mathbb R}^n}\vert\nabla\hat v\vert dx + o(1)\hskip.1cm,
 \end{equation}
 while we easily check that
\begin{equation}\label{Eqt2AfterBUProofTH2}
\int_{{\mathbb R}^n}\vert\nabla\hat v\vert dx 
= \frac{1}{2}\omega_{n-1} + b_{n-1}\hskip.1cm ,
 \end{equation}
where $b_n$ is the volume of the unit ball in ${\mathbb R}^n$. 
Combining (\ref{ProofTH2ContInterEqt})--(\ref{Eqt2AfterBUProofTH2}), we get that the second equation in 
(\ref{TH2Eqt2}) is also true. This ends the proof of Theorem \ref{TH2}.

\medskip The proof of Theorem \ref{TH2}, and hence the theorem itself, 
easily extend to other, more general, types of holes. For instance, when we do not assume anymore that 
only one of the $\varepsilon_{i,\delta}$ is leading, or when we substract 
$K_{i,\delta} \subset B_{x_i}(\varepsilon_{i,\delta})$ instead of the whole 
ball. Only slight modifications in the proof, that we leave to the reader, are required to get such 
extensions.

\section{proof of theorem \ref{TH3}}\label{ProofTH3}

By the change of variables formula for the total variation 
(see Giusti \cite{Giu}), if $T$ is a $C^1$-diffeomorphism from 
${\mathbb R}^n$ to ${\mathbb R}^n$, $\Omega$ is a smooth open subset of ${\mathbb R}^n$, and $u \in BV(\Omega)$, then
\begin{equation}\label{ChangeVarFor}
\int_{\Omega^\star}\vert\nabla u^\star\vert dx 
= \int_\Omega\vert(DT)^{-1}\nu_u\vert \vert DT\vert \vert\nabla u\vert dx\hskip.1cm ,
\end{equation}
where $\Omega^\star = T(\Omega)$, $u^\star = u\circ T^{-1}$, $\nu_u$ is the Radon-Nikodym 
derivative of $\nabla u$ with respect to $\vert\nabla u\vert$, 
and  $\vert DT\vert$ is the absolute value of the 
determinant of $DT$. By (\ref{ChangeVarFor}) with $T = T_\delta$, noting that 
$\vert\nu_u\vert = 1$ for $\vert\nabla u\vert$-almost all $x$, by the variational definition of 
$\lambda_{1,\Omega_\delta}$, and by (\ref{DiffeomDefEqt}), we easily get that
$\limsup_{\delta\to 0}\lambda_{1,\Omega_\delta} \le \lambda_{1,\Omega}$. 
Conversely, we let $u_\delta$ be a nonnegative eigenfunction for $\lambda_{1,\Omega_\delta}$, and we define the 
function $v_\delta$ by 
$v_\delta = u_\delta\circ T_\delta$. By (\ref{DiffeomDefEqt}), (\ref{ChangeVarFor}), and what we just said, 
the sequence $(\overline{v}_\delta)$ is bounded in $BV({\mathbb R}^n)$, with the additional properties that
\begin{equation}\label{ProofTH3Eqt2}
\begin{split}
&\int_{{\mathbb R}^n}\vert\nabla\overline{v}_\delta\vert dx \le \left(1+o(1)\right) \lambda_{1,\Omega_\delta}
\hskip.1cm ,\hskip.1cm\hbox{and}\\
&\int_\Omega v_\delta dx = \left(1 + o(1)\right)\int_{\Omega_\delta}u_\delta dx = 1 + o(1)\hskip.1cm ,
\end{split}
\end{equation}
where $o(1) \to 0$ as $\delta \to 0$. As above, we adopt the notation  $\overline{v}_\delta$ 
for the extension of $v_\delta$ by zero outside $\Omega$. 
Let $D$ be a bounded domain in ${\mathbb R}^n$ which contain both 
$\Omega$ and the $\Omega_\delta$'s. By compactness of the embedding of $BV(D)$ into $L^1(D)$, we may 
assume that, up to a subsequence, $\overline{v}_\delta \to v$ in $L^1(D)$ and almost everywhere as $\delta \to 0$. 
Let $u$ be the restriction of $v$ to $\Omega$. Then $u \ge 0$ and, by the second equation in (\ref{ProofTH3Eqt2}), 
$\int_\Omega udx = 1$. Moreover, by the first equation in (\ref{ProofTH3Eqt2}), and by lower semicontinuity of the 
total variation, 
\begin{eqnarray*}
\lambda_{1,\Omega}
&\le& \int_{{\mathbb R}^n}\vert\nabla\overline{u}\vert dx\\
&\le& \int_D\vert\nabla v\vert dx
\le  \liminf_{\delta \to 0}\lambda_{1,\Omega_\delta}\hskip.1cm ,
\end{eqnarray*}
 where $\overline{u}$ stands for the extension of $u$ by zero outside $\Omega$. 
In particular, $\lambda_{1,\Omega_\delta} \to \lambda_{1,\Omega}$ as $\delta \to 0$, 
and the function $\delta \to \lambda_{1,\Omega_\delta}$ is continuous at $\delta = 0$. This proves the first 
assertion in Theorem \ref{TH3}.

\medskip As a consequence of the above developments, 
$$\int_{{\mathbb R}^n}\vert\nabla\overline{v}_\delta\vert dx 
\to \int_{{\mathbb R}^n}\vert\nabla\overline{u}\vert dx$$ 
as $\delta \to 0$, 
and $u$ is an eigenfunction for $\lambda_{1,\Omega}$. In particular, like in Section \ref{ProofTH1}, 
$\mu_\delta \rightharpoonup \mu$ weakly, where $\mu_\delta = \vert\nabla\overline{v}_\delta\vert$ 
and $\mu = \vert\nabla\overline{u}\vert$. 
In what follows we let $A_\delta$ be an eigenset 
for $\lambda_{1,\Omega_\delta}$ and let $u_\delta = \vert A_\delta\vert^{-1}\chi_{A_\delta}$. Then (we refer again to Section 
\ref{ProofTH1} for the simple argument involved here), $\overline{u} = \vert A\vert^{-1}\chi_A$, where $A$ 
is an eigenset for $\lambda_{1,\Omega}$. In particular, 
$\vert\nabla\overline{v}_\delta\vert \rightharpoonup \mu$ weakly, and 
$\mu = \vert A\vert^{-1} {\mathcal H}^{n-1} \lfloor \partial_\star A$, 
where $\partial_\star A$ is the measure theoretic boundary of $A$ (see for instance Evans-Gariepy \cite{EvaGar}), 
and ${\mathcal H}^{n-1}$ is the $(n-1)$-dimensional Hausdorff measure. Now we prove the second assertion in Theorem \ref{TH3}, 
namely the differentiability of 
the eigenvalue $\lambda_{1,\Omega_\delta}$ at 
$\delta = 0$ and the equation $\left(\lambda_{1,\Omega_\delta}\right)^\prime(0) = \Lambda\lambda_{1,\Omega}$.
As is easily checked from (\ref{DiffeomDefEqt}), 
for any $x \in \overline{\Omega}$, and any $X \in {\mathbb R}^n$ such that 
$\vert X\vert = 1$,
\begin{equation}\label{ValueinChgeVarForm}
\begin{split}
&\vert DT_\delta(x)\vert = 1 - n\Lambda \delta + o(\delta)
\hskip.1cm,\hskip.1cm\hbox{and}\\
& \left\vert\left(DT_\delta\right)(x)^{-1}.X\right\vert = 1 + \Lambda\vert X\vert^2\delta + o(\delta)\hskip.1cm ,
\end{split}
\end{equation}
where the $o(\delta)$'s are uniform in $x$ 
and $X$. By (\ref{SecondDef}) and 
(\ref{ChangeVarFor}) with $T = T_\delta$, we can write that
\begin{equation}\label{ProofTH3FirstMaj}
\lambda_{1,\Omega_\delta} \le 
\frac{\int_{{\mathbb R}^n}\left\vert(DT_\delta)^{-1}.\nu_u\right\vert 
\left\vert DT_\delta\right\vert \vert\nabla\overline{u}\vert dx}
{\int_\Omega\left\vert DT_\delta\right\vert u dx}\hskip.1cm ,
\end{equation}
where $u$ and $\overline{u}$ are as above. By (\ref{ValueinChgeVarForm}),
\begin{equation}\label{ProofTH3DenValue}
\int_\Omega\left\vert DT_\delta\right\vert u dx
= 1 - n\Lambda\delta + o(\delta)\hskip.1cm .
\end{equation}
Since $\vert\nu_u\vert = 1$ a.e w.r.t $\mu$, and ${\mathcal H}^{n-1}(\partial_\star A\backslash\partial^\star A) = 0$, 
we also get with (\ref{ValueinChgeVarForm}) that
\begin{equation}\label{ProofTH3NumValue}
\frac{1}{\lambda_{1,\Omega}} \int_{{\mathbb R}^n}\left\vert(DT_\delta)^{-1}.\nu_u\right\vert 
\left\vert DT_\delta\right\vert \vert\nabla\overline{u}\vert dx
 = 1 - (n-1)\Lambda \delta + o(\delta)\hskip.1cm .
\end{equation}
Plugging (\ref{ProofTH3DenValue}) and (\ref {ProofTH3NumValue}) into (\ref{ProofTH3FirstMaj}), it 
follows that
\begin{equation}\label{ProofTH3DiffFinEqt1}
\lambda_{1,\Omega_\delta} - \lambda_{1,\Omega} 
\le \Lambda\delta\lambda_{1,\Omega} + o(\delta)\hskip.1cm .
\end{equation}
In order to get the converse inequality, we write, still 
using (\ref{SecondDef}) and (\ref{ChangeVarFor}), that
\begin{equation}\label{ProofTH3ConvIneEqt1}
\lambda_{1,\Omega_\delta} - \lambda_{1,\Omega} 
\ge \frac{\int_{{\mathbb R}^n}\left\vert(DT_\delta)^{-1}.\nu_\delta\right\vert 
\left\vert DT_\delta\right\vert \vert\nabla\overline{v}_\delta\vert dx}
{\int_\Omega\left\vert DT_\delta\right\vert v_\delta dx} - 
\frac{\int_{{\mathbb R}^n}\vert\nabla\overline{v}_\delta\vert dx}
{\int_\Omega v_\delta dx}\hskip.1cm ,
\end{equation}
where $\nu_\delta$ is the Radon-Nikodym derivative of $\nabla\overline{v}_\delta$ 
with respect to $\vert\nabla\overline{v}_\delta\vert$. By (\ref{ValueinChgeVarForm}), 
since $\vert\nu_\delta\vert = 1$ for $\mu_\delta$-almost all points,
\begin{equation}\label{PrTH3InterEqtA1}
\int_{{\mathbb R}^n}\left\vert(DT_\delta)^{-1}.\nu_\delta\right\vert 
\left\vert DT_\delta\right\vert \vert\nabla\overline{v}_\delta\vert dx
= \int_{{\mathbb R}^n}\vert\nabla\overline{v}_\delta\vert dx 
- (n-1)\Lambda\delta \int_{{\mathbb R}^n}\vert\nabla\overline{v}_\delta\vert dx 
+ o(\delta)
\end{equation}
and
\begin{equation}\label{PrTH3InterEqtA2}
\int_\Omega\left\vert DT_\delta\right\vert v_\delta dx = \int_\Omega v_\delta dx 
- n \Lambda \delta \int_\Omega v_\delta dx + o(\delta)\hskip.1 cm.
\end{equation}
Combining (\ref{ProofTH3ConvIneEqt1}), (\ref{PrTH3InterEqtA1}), and (\ref{PrTH3InterEqtA2}), it follows that
\begin{equation}\label{ProofTH3RevIneqEqtSte1}
\lambda_{1,\Omega_\delta} - \lambda_{1,\Omega} \ge 
\frac{\int_{{\mathbb R}^n}\vert\nabla\overline{v}_\delta\vert dx}{\int_\Omega v_\delta dx} 
\Lambda \delta + o(\delta)\hskip.1cm .
\end{equation}
Since  $\int_\Omega v_\delta dx \to 1$ and 
$\int_{{\mathbb R}^n}\vert\nabla\overline{v}_\delta\vert dx 
\to \lambda_{1,\Omega}$ as $\delta \to 0$, it follows from 
(\ref{ProofTH3DiffFinEqt1}) and (\ref{ProofTH3RevIneqEqtSte1}) that
$\lambda_{1,\Omega_\delta} - \lambda_{1,\Omega} = 
\Lambda \delta \lambda_{1,\Omega} + o(\delta)$.
The equation holds for a subsequence, but since the right hand side in 
the equation does not depend on the subsequence, it holds true for all $\delta$. 
In particular, $\left(\lambda_{1,\Omega_\delta}\right)^\prime(0) = \Lambda\lambda_{1,\Omega}$ and this 
ends the proof of Theorem \ref{TH3}.

\medskip The proof of the first assertion in Theorem 
\ref{TH3}, and hence the continuity of 
$\lambda_{1,\Omega_\delta}$ at $\delta = 0$, extend to 
very general $T_\delta$'s. We basically only need that $T_\delta \to Id$ in 
the $C^1$-topology as $\delta \to 0$.

\end{document}